\documentclass[a4paper,11pt]{amsart}

\usepackage{a4wide}
\usepackage{amsmath,amssymb,amsthm}
\usepackage{amsfonts,amsthm,amsmath,amssymb,amscd,color}
\usepackage{enumerate}
\usepackage{parskip}
\usepackage{graphicx}
\usepackage{babel}

\newcommand{\RR}{\mathbb{R}}

\newcommand{\NN}{\mathbb{N}}

\newcommand{\eps}{\varepsilon}

\newcommand{\CB}{{\mathcal{B}}}

\newcommand{\CP}{\mathcal{P}}
\newcommand{\CM}{{\mathcal{M}}}

\newcommand{\BL}{\mathrm{BL}}
\newcommand{\TV}{\mathrm{TV}}
\newcommand{\FM}{\mathrm{FM}}

\newcommand{\supp}{\mathrm{supp}}

\newtheorem{theo}{Theorem}[section]
\newtheorem{prop}{Proposition}[section]
\newtheorem{lemma}{Lemma}[section]

\newtheorem{corollary}{Corollary}[section]

\newtheorem{example}{Example}[section]
\newtheorem{remark}{Remark}[section]

\newcommand{\smfrac}[2]{\mbox{$\frac{#1}{#2}$}}

\newcommand{\mathbbm}[1]{{#1\!\!#1}}
\newcommand{\ind}{{\mathbbm{1}}}
\newcommand{\Mol}{{\mathcal{M}ol}}

\begin{document}
	\title[Fortet-Mourier distance to weighted sums of Dirac measures]{Explicit expressions and computational methods for the Fortet-Mourier distance to finite weighted sums of Dirac measures}
	\author{Sander C. Hille}
	\address{Mathematical Institute, Leiden University, P.O. Box 9512, 2300 RA Leiden, The Netherlands, (SH)}
	\email{shille@math.leidenuniv.nl (corresponding author)}
	
	\author{Esm\'ee S. Theewis}
	\address{Mathematical Institute, Leiden University, P.O. Box 9512, 2300 RA Leiden, The Netherlands, (ET)}
	\email{e.s.theewis@umail.leidenuniv.nl}
	
	\date\today
	
	\subjclass[2020]{28A33, 46E27, 90C05, 90C25} % 92B05,37A30,60J20
	\keywords{Fortet-Mourier norm, Borel measure, metric space, linear and convex optimization}
	
	\begin{abstract}
		Explicit expressions and computational approaches are given for the Fortet-Mourier distance between a positively weighted sum of Dirac measures on a metric space and a positive finite Borel measure. Explicit expressions are given for the distance to a single Dirac measure. For the case of a sum of several Dirac measures one needs to resort to a computational approach. In particular, two algorithms are given to compute the Fortet-Mourier norm of a molecular measure, i.e. a finite weighted sum of Dirac measures. It is discussed how one of these can be modified to allow computation of the dual bounded Lipschitz (or Dudley) norm of such measures.
	\end{abstract}
	
	\maketitle

\section{Introduction}

Let $(S,d)$ be a metric space, equipped with its Borel $\sigma$-algebra $\CB(S)$.  We denote by $\BL(S)$ the real vector space of bounded Lipschitz functions on $(S,d)$. The Lipschitz constant of $f\in\BL(S)$ is written as $|f|_L$. Following Lasota, Szarek and co-workers (e.g. \cite{Lasota-Myjak, Lasota-Myjak-Szarek}) we define the {\it Fortet-Mourier norm} on the finite signed Borel measures $\CM(S)$ on $S$ by
\begin{equation}\label{eq:definition FM-norm}
	\|\mu\|_\FM^* := \sup\bigl\{ \langle \mu,f\rangle: f\in \BL(S),\ \|f\|_\infty\leq 1,\ |f|_L\leq 1\bigr\},\qquad \mu\in\CM(S),
\end{equation} 
where the indicated pairing is given by integration: $\langle\mu,f\rangle := \int_S f\,d\mu$. In this paper we provide explicit expressions and computational methods for Fortet-Mourier norms of the form 
\begin{equation}\label{eq:FM-norms of interest}
	\bigl\|\sum_{i=1}^N \alpha_i\delta_{x_i}-\mu\bigl\|_\FM^*,\qquad N\in\NN,\ \alpha_i>0,\ x_i\in S\ \mbox{and}\ \mu\in \CM^+(S).
\end{equation}
Here, $\delta_x$ denotes the Dirac (or point) measure located at $x\in S$ and $\CM^+(S)$ is the convex cone of positive measures in $\CM(S)$. 

This norm, or the equivalent dual bounded Lipschitz norm (also called Dudley norm or flat metric -- for the derived metric), is used much in the study of dynamical systems in spaces of measures.
For example, one encounters these norms in the context of Markov operators and semigroups on (probability) measures \cite{Czapla_ea:2019,Hille-Worm:2009,Lasota-Myjak-Szarek}, like those defined by Iterated Function Systems \cite{Lasota-Myjak} or Piecewise Deterministic Markov Processes \cite{Alkurdi-Hille-vanGaans:2013,Czapla_ea:2019, Hille-Horbacz-Szarek}. Deterministic systems in spaces of measures appeared e.g. in models for population dynamics and biological systems \cite{Azmy_ea:1999,Azmy_ea:2005,Carillo-Colombo-Gwiazda:2012}, transport equations \cite{Ackley-Saintier_ea:2019,Gwiazda-Hille_ea:2019,Piccoli-Rossi-Tournus:2021} and interacting particle systems or crowd dynamics \cite{Evers-H-M, Piccoli-Tosin:2011}.

In a setting where the dynamics conserve total mass, many authors have used the much more studied family of Wasserstein distances \cite{Villani:2003, Piccoli-Rossi:2014}. These distances are defined for measures of positive equal mass only, however. If total mass can vary, as in most of the mentioned examples of deterministic type, then Wasserstein metrics are of limited use. Extensions are being explored \cite{Piccoli-Rossi-Tournus:2021}, but Fortet-Mourier, Dudley or flat metric -- or a metric defined by duality with a H\"older class of functions \cite{Gwiazda-Hille_ea:2019} -- may be preferred.

In settings with varying total mass, norms of the form \eqref{eq:FM-norms of interest} are of interest for several reasons. In a measure framework, continuum models and discrete interacting particle descriptions can be framed within one functional analytic setting. A weighted sum of Dirac measures then represents the particle model, while a measure $\mu$ that is absolutely continuous with respect to Lebesgue measure represents the other. Norm \eqref{eq:FM-norms of interest} then quantifies the deviation between these two descriptions. For example, in so-called Patlak-Keller-Segel type chemotactic models it has been shown that the continuum solution converges to sums of Dirac measures in finite time \cite{Herrero-Velazquez:1997, Dolbeault-Schmeiser:2009}, yielding blow-up in the used $L^p$-norm ($p>1$). Expressions like \eqref{eq:FM-norms of interest} may trace such `concentration of mass' and express a rate of convergence.\\
In numerical analysis of particular continuum models it may be advantageous to simulate a well-chosen interacting particle system instead of simulating the partial differential equations. See e.g. \cite{Cottet-Mas-Gallic:1990} where this was advocated for the two-dimensional Navier-Stokes problem, to limit numerical diffusion. Estimates of norms of the form \eqref{eq:FM-norms of interest} then appear naturally in error estimates.\\
Within the single setting of a particle model, a question is to quantify deviation between two instances of the model, with different particle number. Then, $\mu$ is also a weighted sum of Dirac measures. Expression \eqref{eq:FM-norms of interest} then reduces to computing norms of the form
\begin{equation}\label{eq:FM-norm of molecular measure}
	\|\tau\|_\FM^*,\qquad \mbox{with}\ \tau\in\Mol(S):=\mathrm{span}_{\RR}\bigl\{ \delta_x:x\in S\bigr\},
\end{equation}
the subspace of so-called {\it molecular measures} \cite{Pachl}.

There exist a few results that provide exact algorithms to compute norms of molecular measures. Jab{\l}o\'nski and Marciniak-Czochra \cite{Jablonski-MCz:2013} provided an algorithm to compute $\|\tau\|_\FM^*$ with $\tau\in\Mol(S)$ and $S=\RR$ with the Euclidean metric or a bounded closed interval therein (see also \cite{Evers:thesis} Appendix, for a description and application of their algorithm). Their approach depends heavily on the total ordering that is available on $\RR$. Generalization of this approach to higher dimension or to any Polish space is therefore inhibited. Sriperumbudur {\it et al.} \cite{Sriperumbudur_ea:2012} provided an algorithm for computing the (equivalent) Dudley norm of a difference of two empirical measures. That is,  $\tau=\nu-\mu$ with $\nu=\smfrac{1}{N}\sum_{i=1}^N \delta_{x_i}$ and a similar expression for $\mu$, possibly with a different number of point measures (see \cite{Sriperumbudur_ea:2012} Theorem 2.3, p.1557). The state space $S$ can be any metric space. 

Up till now, to our knowledge, neither for a specific choice of $S$, nor in the generality of an arbitrary metric space $(S,d)$, there are hardly any explicit expressions for \eqref{eq:FM-norms of interest}, except for the well-known
\begin{equation}\label{eq:well-known norms}
	\bigl\|\delta_x-\delta_y \bigr\|^*_\FM = 2\wedge d(x,y),\qquad x,y\in S
\end{equation}
(see e.g. \cite{Hille-Worm:2009, Pachl}). Our main results, Theorem \ref{prop: general distance with thetas} and Theorem \ref{prop:FM-norm Dirac minus positive}, will allow to extend this e.g. to the expression
\begin{equation}\label{eq:distance Dirac to pos measure}
	\bigl\| \delta_x - \mu\bigr\|_\FM^* = \bigl\langle \mu, 2\wedge d(x,\cdot)\bigr\rangle,\qquad x\in S,\ \mu\in \CP(S),
\end{equation}
where $\CP(S)$ is the subset of probability measures in $\CM^+(S)$ (see Proposition \ref{clry: FM-norm Dirac min prob} and various corollaries of Theorem \ref{prop:FM-norm Dirac minus positive} in Section \ref{sec:FM-distance to single Dirac}), or the expression
\begin{equation}
	\bigl\| \alpha\delta_x - \beta\delta_y \bigr\|_\FM^* = |\alpha-\beta|\ +\ (\alpha\wedge\beta)(2\wedge d(x,y), \qquad \alpha,\beta>0,\ x,y\in S.
\end{equation}
(See Corollary \ref{clry:FM-distance two weighted single Diracs}). Such explicit expressions may be useful in obtaining (better) estimates of Fortet-Mourier norms of the indicated form. Moreover, expression \eqref{eq:distance Dirac to pos measure} enables the explicit computation of e.g. the Fortet-Mourier distance of a Dirac measure to a measure that is absolutely continuous with respect to Lebesgue measure, which is a novel result.

Another motivation for determining expressions for norms like \eqref{eq:FM-norms of interest} comes from approximation theory. The mathematical question in which these explicit formulae may be of help is in that of existence and computation of {\it best approximation} of $\mu\in\CM^+(S)$ by a a positive sum of at most $N$ Dirac measures in Fortet-Mourier distance, where $N$ is fixed {\it a priori}. This is e.g. relevant for an interacting particle approach to solving a continuum model. The continuum initial condition must then be replaced by a number of particles. How can these be `best' distributed over space, such that the error caused by the approximation of the initial condition is minimal? Is there such a best approximation? Can it be found computationally?

General results for the existence of a best approximation have been known for long, e.g. for reflexive Banach spaces and closed convex sets therein, see e.g. \cite{Deutsch:1980}. Although the indicated set of sums of $N$ Dirac measures is closed, it is not convex. Moreover, the completion of the space $\CM(S)$ for the $\|\cdot\|_\FM^*$-norm is hardly ever reflexive. (In fact, it is isometrically isomorphic to $\BL(S)$ with norm $\|\cdot\|_\FM:=\max(\|f\|_\infty,|f|_L)$, which can be proven in similar way as \cite{Hille-Worm:2009}, Theorem 3.7, p.360). Nevertheless, a best approximation can be shown to exist on compact and complete metric spaces, essentially by exploiting the compactness of the space of probability measures that is provided by Prokhorov's Theorem.

On non-compact spaces, the situation is much more delicate, as can be illustrated by the following particular case. Expression \eqref{eq:distance Dirac to pos measure} allows to reformulate the special case $N=1$ and $\mu$ a {\it convex} sum of finitely many Dirac measures located at $x_j$ to the problem of minimizing over $x\in S$ the expression
\begin{equation}
	\bigl\|\delta_x - \sum_{j=1}^n \alpha_j\delta_{x_j}\bigr\|_\FM^* = \sum_{j=1}^n \alpha_j\, (2\wedge d(x,x_j)),\qquad (\alpha_j>0,\ \mbox{$\sum_j$}\alpha_j = 1).
\end{equation}
This minimization problem is the Fermat-Weber problem for the metric $d':=2\wedge d$ on $S$, with weights $\alpha_j$. In economic terms, a solution to Weber's problem provides an optimal location $x$ for a production site such that products produced there, can be distributed to the distribution sites at $x_j$ with minimal cost, when the $\alpha_j$ represent the transport cost per item per unit distance \cite{Drezner_ea:2002}. Fermat's original problem was to construct a point $x$ for which the sum of the distances to three given points is minimal. For $n>3$ and equal weights, in Euclidean space, there does not exist a geometric construction for the best point. Existence of a minimizer to the general Fermat-Weber problem is guaranteed for so-called Hadamard spaces (also called complete $\mathrm{CAT}(0)$ spaces) by \cite{Bacak:2014}, Lemma 2.2.19. Various numerical schemes have been developed to determine a minimizer, e.g. \cite{Kuhn-Kuenne:1962,Drezner_ea:2002}. Research on the Fermat-Weber problem continues to this date \cite{Church-Drezner_ea:2022}.

In Section \ref{sec:dimensional reduction} and Section \ref{sec:FM-distance to single Dirac} we present our main results on explicit expressions, Theorem \ref{prop: general distance with thetas} and Theorem \ref{prop:FM-norm Dirac minus positive}, and various consequences derived from these. We shall present an algorithmic approach to computing Fortet-Mourier norms of the form \eqref{eq:FM-norms of interest} with $\nu$ a positive molecular measure in Section \ref{sec:algorithms}. Section \ref{sec:computing FM-norm molecular} is concerned with algorithms for computing $\|\tau\|_\FM^*$ for any $\tau\in\Mol(S)$. In both sections, $S$ is assumed to be a metric space, without additional constraints, like separability or completeness. This substantially generalizes \cite{Jablonski-MCz:2013}. Section \ref{sec:results for Dudly norm} discusses how the results of Section \ref{sec:computing FM-norm molecular} can be modified to compute the Dudley norm of any molecular measure. This generalizes the result in \cite{Sriperumbudur_ea:2012} on this topic.

\subsection{Preliminary results and notation}

For a metric space $(S,d)$ we let $\BL(S)$ denote the ordered vector space of real-valued bounded Lipschitz functions on $S$, with point-wise partial order. We suppress the metric in notation, because there will be no need to consider multiple metrics on the same space. For $f\in\BL(S)$, 
\[
|f|_L := \sup\bigl\{ \frac{|f(x)-f(y)|}{d(x,y)}: x,y\in S,\ x\neq y\bigr\}
\]
denotes the Lipschitz constant of $f$. Occasionally, we shall write $|f|_{L,S}$ if we wish to stress the underlying metric space. $\BL(S)$ is an algebra for point-wise multiplication. It is also an ordered vector space for point-wise ordering, even a vector lattice (Riesz space) with supremum $f\vee g$ and infimum $f\wedge g$ of two elements $f,g\in\BL(S)$ given by point-wise maximum and minimum:
\[
(f\vee g)(x) = \max\bigl(f(x),g(x)),\qquad (f\wedge g)(x) =\min\bigl( f(x),g(x)\bigr),\qquad x\in S.
\]
One has 
\[
|f\vee g|_L\leq \max(|f|_L,|g|_L)\quad \mbox{and}\quad |f\wedge g|_L\leq \max(|f|_L,|g|_L),
\]
see e.g. \cite{Dudley:1966}.
We consider the norm on $\|f\|_\FM := \max\bigl(\|f\|_\infty,|f|_L)$ on $\BL(S)$, which turns $\BL(S)$ into a Banach space. The unit ball in $\BL(S)$ for this norm is denoted by $B^S_\FM$.

The following lemma yields a result on a recurring construction related to $B^S_\FM$ that will appear in several proofs later. \begin{lemma}\label{lem:bounding from below by extreme point} 
	Let $(S,d)$ be a metric space. The following statemenst hold:
	\begin{enumerate}
		\item[{\it (i)}] Let $g\in B^S_\FM$, $N\in \NN$ and let $x_i\in S$ ($i=1,\dots, N$). Define
		\begin{equation}\label{eq:def h}
			h := (-\ind)\vee\bigvee_{i=1}^N (g(x_i) - d(x_i,\cdot)).
		\end{equation}
		Then $h\in B^S_\FM$, $h\leq g$ and $h(x_i)=g(x_i)$ for all $i\in\{1,\dots,N\}$.
		
		\item[{\it (ii)}] Let $P=\{x_1,\dots,x_N\}$ be a subset of $S$ of distinct points, equipped with the restriction of the metric on $S$. If $g\in B^P_\FM$, then $h$ defined by \eqref{eq:def h} is in $B^S_\FM$ and satisfies $h(x_i)=g(x_i)$ for all $i$. Moreover,  $|h|_{L,S} \leq |g|_{L,P}$.
	\end{enumerate}
	
\end{lemma}
\begin{proof}
	{\it (i)}.\ The functions, $x\mapsto g(x_i)-d(x_i,x)$ are Lipschitz on $S$ with Lipschitz constant at most 1 and bounded from above by 1, since $\|g\|_\infty\leq 1$. Hence $-\ind\leq h\leq \ind$ and $|h|_L\leq 1$. Thus, $h\in B^S_\FM$.
	First we show that $h\leq g$. Take $x\in S$. Then either $h(x)=-1$ or $h(x)=g(x_i)-d(x_i,x)$ for some $i\in\{1,\dots,N\}$. In the first case, one trivially has $h(x)=-1\leq -\|g\|_\infty\leq g(x)$. In the other case, one has $g(x_i)-g(x)\leq |g|_Ld(x_i,x) \leq d(x_i,x)$. So $h(x)=g(x_i)-d(x_i,x)\leq g(x)$.
	Next, by construction of $h$ it holds for every $i\in\{1,\dots,N\}$ that $h(x_i)\geq g(x_i) - d(x_i,x_i)= g(x_i)$. Thus $h(x_i)=g(x_i)$.
	
	\noindent {\it (ii)}.\ Follows from part {\it (i)} and the McShane Extension Theorem, \cite{McShane} Theorem 1.
\end{proof}

The following lemma gives two elementary properties of the maximum operator, which will be useful in Section \ref{sec:algorithms}. 
\begin{lemma}\label{lem:properties max operator}
	Let $n\in\NN$ and $a_i,b_i\in\RR$ for $i=1,\dots,n$. Then
	\begin{enumerate}
		\item[{\it (i)}] $\max_i (a_i+b_i) \leq \max_i a_i + \max_i b_i$,
		\item[{\it (ii)}] $\bigl| \max_i a_i - \max_i b_i\bigr| \leq \max_i |a_i-b_i|$.
	\end{enumerate}
\end{lemma}
\begin{proof}
	Part {\it (i)} is obvious. For part {\it (ii)}, without loss of generality, assume that $\max_i a_i\geq \max_i b_i$. let $j\in\{1,\dots,n\}$ be such that $a_j = \max_i a_i$. Then
	\[
	\bigl| \max_i a_i - \max_i b_i\bigr|\ =\ a_j - \max_i b_i\ \leq\  a_j - b_j\ \leq\ |a_j-b_j|\ \leq\ \max_i |a_i-b_i|.
	\]
\end{proof}

For a metric space $(S,d)$, $\CM(S)$ embeds naturally into the dual space $\BL(S)^*$ of continuous linear functionals on $\BL(S)$ by means of the `integration functional with respect to $\mu$', $I_\mu$:
\[
I_\mu(f) := \int_S f\,d\mu =: \langle \mu, f\rangle,\qquad \mu\in\CM(S),\ f\in\BL(S). 
\]
The map $\mu\mapsto I_\mu$ is injective, because the indicator function of any closed set $C\subset S$ can be approximated point-wise by a decreasing sequence of functions in $\BL(S)$, namely $f_n:= [1-nd(\cdot, C)]^+$, and any $\mu\in\CM(S)$ is regular, because $S$ is a metric space (cf. \cite{Bogachev-II:2007} Theorem 7.17). By means of the mentioned embedding, the norm $\|\cdot\|_\FM$ introduce a norm on $\CM(S)$ through the dual space $\BL(S)^*$, which is precisely given by \eqref{eq:definition FM-norm} and which can be found in part of the literature under the name `Fortet-Mourier norm'.

%The support of $\mu\in\CM(S)$ is
%\[
%\supp(\mu) := \bigcap \bigl\{ C\subset S: C\ \mbox{is closed},\ |\mu|(C) = |\mu|(S)\bigl\}.
%\]
%In general, $|\mu|\bigl(\supp(\mu)\bigr)\leq |\mu|(S)$. One says that `{\it $\mu$ has support}' if $|\mu|$ is concentrated on $\supp(\mu)$, i.e. equality holds. If $S$ is a separable metric space, then any $\mu\in\CM(S)$ has support (cf. \cite{Bogachev-II:2007} Proposition 7.2.9, p.77). Then one can show that
%\begin{equation}\label{eq:alternative description support}
%\supp(\mu) = \bigl\{ x\in S: |\mu|(B(x,r))>0\ \mbox{for all}\ r>0 \bigr\},
%\end{equation}
%making use of the Lindel\"of property of a separable metric space. Expression \eqref{eq:alternative description support} has been used as definition of support for measures on Polish spaces (cf. eg. \cite{Lasota-Myjak}). On a general metric space, separable measures in $\CM(S)$, i.e. those that are concentrated on a separable set,  have support.

%We denote by $\CM^+(S)$ the convex cone of positive measures in $\CM(S)$ and by $\CP(S)$ the convex set of probability measures. The Dirac measure or point-mass located at $x\in S$ is denoted by $\delta_x$. A finite linear combination of Dirac measures is called a {\it molecular measure}. The vector space of molecular measures is written as $\Mol(S)$, while $\Mol^+(S)$ consists of the positive measures. 

\section{Dimensional reduction for determining the defining supremum}
\label{sec:dimensional reduction}

Let $(S,d)$ be a metric space. Neither completeness, nor separability is required. The defining expression for the Fortet-Mourier norm \eqref{eq:definition FM-norm} cannot be conveniently used for the computation of this norm in practice. The main issue is, that there is no convenient method (yet) to determine the supremum over the full unit ball $B^S_\FM$.

The following key result allows to substantially reduce the dimension of the set over which to take the supremum, provided one of the two measures is molecular.
%It has a spirit similar to Proposition \ref{prop: FM-distance probability measures}:
\begin{theo}\label{prop: general distance with thetas}
	Let $\nu\in\Mol^+(S)$ with $P:=\supp(\nu)=\{x_1,\dots,x_N\}$, with the $x_i$ all distinct. Let $\mu\in\CM^+(S)$. Then
	\begin{align}
		\bigl\|\nu-\mu\bigr\|^*_\FM\ & =\ \sup_{f\in B^P_\FM} \bigl \langle\nu-\mu, (-\ind)\vee \bigvee_{i=1}^N (f(x_i)- d(x_i,\cdot)) \bigr\rangle\label{eq:FM-norm via BP}\\
		& = \ \sup_{\theta\in[-1,1]^N}  \bigl\langle\nu-\mu, (-\ind)\vee \bigvee_{i=1}^N (\theta_i- d(x_i,\cdot)) \bigr\rangle.\label{eq:FM-norm via thetas}
	\end{align}
\end{theo}
\begin{proof}
	First of all, $\|\nu-\mu\|_\FM^* = \sup_{g\in B^S_\FM} \langle\nu-\mu,g\rangle$. Moreover, $\nu=\sum_{i=1}^N \alpha_i \delta_{x_i}$, with $\alpha_i>0$. Let $g\in B^S_\FM$. Note that its restriction to $P$, $g|_P$, is in $B^P_\FM$ and $g(x_i)=g|_P(x_i)$ for all $i\in\{1,\dots,N\}$. Define  
	\begin{equation}\label{def:h}
		h := (-\ind)\vee \bigvee_{i=1}^N (g(x_i)-d(x_i,\cdot)).
	\end{equation}
	According to Lemma \ref{lem:bounding from below by extreme point}, $h\in B^S_\FM$, $h\leq g$ and $h(x_i)=g(x_i)$ for all $i$. Therefore,
	\begin{equation}\label{eq:inequality measures on h and g}
		\langle\nu-\mu,h\rangle\ =\ \sum_{i=1}^N \alpha_i g(x_i) - \langle\mu,h\rangle\ \geq \ \sum_{i=1}^N \alpha_i g(x_i) - \langle\mu,g\rangle\ =\ \langle\nu-\mu,g\rangle.
	\end{equation}
	This proves inequality `$\leq$' in \eqref{eq:FM-norm via BP}. The other inequality in this equation is an immediate consequence of the observation that for any $f\in B^P_\FM$ the function $(-\ind)\vee\bigvee_{i=1}^N (f(x_i) - d(x_i,\cdot))$ is in $B^S_\FM$ (using Lemma \ref{lem:bounding from below by extreme point} {\it (ii)}).
	
	For proving equality \eqref{eq:FM-norm via thetas}, let $\theta\in[-1,1]^N$ and define
	\begin{equation}\label{def:h_theta}
		h_\theta (x):= (-\ind)\vee \bigvee_{i=1}^N (\theta_i-d(x_i,x)),\qquad x\in S.
	\end{equation}
	By definition, $-1\leq h_\theta\leq 1$. Moreover, $|h_\theta|_L \leq \max_{1\leq i\leq N}\bigl(|-\ind|_L, |\theta_i-d(x_i,\cdot)|_L\bigr)= 1$. Thus, $h_\theta\in B^S_\FM$. So `$\geq$' holds in \eqref{eq:FM-norm via thetas}. 
	For `$\leq$': let $f\in B^P_\FM$ and define
	\[
	g(x) := (-\ind)\vee \bigvee_{i=1}^N (f(x_i)- d(x_i,x))
	\]
	According to Lemma \ref{lem:bounding from below by extreme point} {\it (ii)}, $g\in B^S_\FM$ and $g(x_i)=f(x_i)$ for all $i$. Thus, if $\theta_i:=f(x_i)$, $i=1,\dots, N$, then $\theta_i\in[-1,1]$ and $h_\theta = g$. So, the supremum in \eqref{eq:FM-norm via thetas} is over a larger set than that in \eqref{eq:FM-norm via BP}. Hence, `$\leq$' holds.  
\end{proof}

The difference between the similarly looking functions in \eqref{eq:FM-norm via BP} and \eqref{eq:FM-norm via thetas} is, that $h_\theta(x_i)\geq \theta_i$ for all $i$ in \eqref{eq:FM-norm via thetas}, but equality need not hold, while for $g\in B^P_\FM$, the function $h$ defined by \eqref{def:h}, which is used in \eqref{eq:FM-norm via BP}, does satisfy $h(x_i)=g(x_i)$ for all $i$.
\vskip 0.1cm

If we restrict our attention to $\nu$ and $\mu$ being probability measures, the dimension of the set over which one takes the supremum in \eqref{eq:FM-norm via thetas} can be further reduced by one, as the following result ascertains. 
\begin{prop}\label{prop: FM-distance probability measures}
	Let $\nu\in\Mol^+(S)\cap\CP(S)$, with  $\supp(\nu) = \{x_1,\dots, x_N\}$, all $x_i\in S$ distinct, and let $\mu\in\CP(S)$. Then
	\begin{equation}\label{eq:expresion norm probability meaures}
		\|\nu-\mu\|_\FM^* = \sup\left\{\ \bigl\langle \nu-\mu, (-\ind) \vee \bigvee_{i=1}^N (\theta_i - d(x_i,\cdot)) \bigr\rangle :\  \theta_1,\dots,\theta_N\in[-1,1],\ \max_i(\theta_i)=1\ \right\}.
	\end{equation}
\end{prop}
\begin{proof}
	The inequality `$\geq$' follows immediately from \eqref{eq:FM-norm via thetas}, since the supremum in \eqref{eq:expresion norm probability meaures} is taken over a subset of that in \eqref{eq:FM-norm via thetas}.
	For the other inequality, let $g\in B^S_\FM$. Put $\eps:=\min_{1\leq i\leq N} (1-g(x_i))\geq 0$ and $\theta_i:= g(x_i)+\eps$ for $i=1,\dots,N$. Note that $\theta_i\in[-1,1]$. Let $k\in\{1,\dots,N\}$ be such that $\eps = 1-g(x_k)$. Then $\theta_k=1$, so $\max_i(\theta_i) = 1$. With $\theta=(\theta_1,\dots, \theta_N)$ and $h_\theta$ as defined in \eqref{def:h_theta},
	\begin{equation}\label{eq:estimate h theta below}
		h_\theta(x_i) \geq \theta_i-d(x_i,x_i) = g(x_i)+\eps.
	\end{equation}
	We claim that $h_\theta\leq g+\eps$. Indeed, if $x\in S$ is such that $h_\theta(x)=-1$, then the inequality holds trivially, since $g\geq -1$ and $\eps\geq0$. If $h_\theta(x)>-1$, then for some $j\in\{1,\dots,N\}$,
	\begin{align}
		h_\theta(x)\ &=\ \theta_j - d(x_j,x)\ =\ g(x_j) + \eps - d(x_j,x)\ =\ g(x) + \eps - d(x_j,x) + g(x_j) - g(x)\nonumber\\
		& \leq\ g(x) + \eps - d(x_j,x) + d(x_j,x)\ =\ g(x)+\eps, \label{eq:estimate h theta above}
	\end{align}
	because $g\in B^S_\FM$ and consequently, $|g|_L\leq 1$. Write $\nu=\sum_{i=1}^N \alpha_i\delta_{x_i}$, with $\alpha_i>0$, $\sum_{i=1}^N \alpha_i=1$. The two bounds \eqref{eq:estimate h theta below} and \eqref{eq:estimate h theta above} yield
	\begin{align*}
		\langle \nu-\mu, h_\theta\rangle\ &=\ \sum_{i=1}^N \alpha_i h(x_i) - \int_S h_\theta\,d\mu\ \geq\ \sum_{i=1}^N \alpha_i \bigl( g(x_i)+\eps \bigr) - \int_S (g+\eps)\,d\mu\\
		& =\ \sum_{i=1}^N \alpha_i g(x_i) +\eps \sum_{i=1}^N \alpha_i - \int_S g\,d\mu - \eps\mu(S) \ = \ \langle \nu-\mu, g\rangle.
	\end{align*}
	Since $\|\nu-\mu\|_\FM^* = \sup_{g\in B^S_\FM} \langle \nu-\mu, g\rangle$, we obtain inequality `$\leq$' in \eqref{eq:expresion norm probability meaures}.
\end{proof}

\noindent These results give rise to novel explicit expressions for the Fortet-Mourier distance to a single point mass.

\section{Explicit expressions for the distance to a single point mass}
\label{sec:FM-distance to single Dirac}

Theorem \ref{prop: general distance with thetas} and Proposition \ref{prop: FM-distance probability measures} reduce the supremum expression for the distance to a positive molecular measure to a maximization problem of a suitable continuous function over a particular compact set. In this section we show, that if $\nu$ is a single (weighted) Dirac measure, a location where this maximum is attained can be explicitly determined. Moreover, it will become clear, that this location need not be unique. It results into various novel explicit expressions for norms of the form $\|\alpha\delta_x - \mu\|_\FM^*$, for $x\in S$, $\alpha>0$ and $\mu\in\CM^+(S)$.
\vskip 0.1cm

\begin{prop}\label{clry: FM-norm Dirac min prob}
	Let $x\in S$ and $\mu\in\CP(S)$. Then
	\[
	\|\delta_x-\mu\|^*_\FM\ =\ \bigl\langle \delta_x-\mu, (-\ind)\vee(1-d(x,\cdot)) \bigr\rangle\ =\ \bigl\langle\mu, 2\wedge d(x,\cdot)\bigr\rangle.
	\]
\end{prop}
\begin{proof}
	Specifying the result of Proposition \ref{prop: FM-distance probability measures} to the case $N=1$ yields the first equality. Then observe that
	\begin{equation}\label{eq:rewriting function}
		(-\ind)\vee(1-d(x,\cdot)) = -\bigl( \ind \wedge (d(x,\cdot) - 1)\bigr) = -\bigl( 2\wedge d(x,\cdot) \bigr) + \ind.
	\end{equation}
	Applying the measure $\delta_x-\mu$ to the latter function gives the result, since $\mu\in\CP(S)$.
\end{proof}

\begin{example}
	1.) Applying Proposition \ref{clry: FM-norm Dirac min prob} in the special case $\mu=\delta_y$ yields the well-known expression $\|\delta_x-\delta_y\|^*_\FM=2\wedge d(x,y)$, shown in \eqref{eq:well-known norms}.
	\vskip 0.1cm
	
	\noindent 2.) If $\mu=\alpha\delta_{y_1} + (1-\alpha)\delta_{y_2}$ with $y_1, y_2\in S$ and $0\leq \alpha\leq 1$, then Proposition \ref{clry: FM-norm Dirac min prob} gives the explicit expression
	\begin{equation}\label{eq:FM-norm 1-2 Diracs}
		\bigl\| \delta_x - \alpha\delta_{y_1} - (1-\alpha)\delta_{y_2}\bigr\|^*_\FM\ =\ \alpha(2\wedge d(x,y_1))\ +\ (1-\alpha)(2\wedge d(x,y_2)).
	\end{equation}
	\vskip 0.1cm
	
	\noindent 3.) Let $S=[0,1]$, equipped with the Euclidean metric and let $\lambda$ be the Borel-Lebesgue measure on $S$, $\lambda([0,1])=1$. Then for any $x\in S$,
	\begin{equation}\label{eq:distance Dirac Lebesgue interval}
		\bigl\| \delta_x - \lambda\bigr\|_\FM^*\ =\ \int_0^1 2\wedge |x-y|\,dy\ =\ \int_0^x (x-y) dy + \int_x^1 (y-x)dy\ = \ \mbox{$\frac{1}{2}$} - x + x^2.
	\end{equation}
	Notice that expression \eqref{eq:distance Dirac Lebesgue interval} is minimal for $x=1/2$ with value $1/4$. Thus, there exists a (unique) best approximation in $\CP(S)$ of $\lambda$ by a single Dirac measure in Fortet-Mourier:
	\begin{equation}
		\inf_{x\in S} \bigl\| \delta_x - \lambda\bigr\|_\FM^*\ =\ \bigl\| \delta_{1/2} - \lambda\bigr\|_\FM^*\ =\ \frac{1}{4}.
	\end{equation}
	The location of best approximation is in this case the median of the uniform distribution on $[0,1]$, which is given by $\lambda$.
	\vskip 0.1cm
	
	\noindent 4.) Let $S$ be as in part 3.) and $f\geq 0$ a probability distribution function. Then
	\[
	\bigl\|\delta_x -f\,d\lambda \bigr\|_\FM^* = \int_0^1 |x-y| f(y)dy = \int_0^x (x-y)f(y)dy - \int_x^1 (y-x)f(y)dy.
	\]
	For specific $f$ the latter expression can be computed, in principle.
\end{example}

The following result is an immediate corollary of Theorem \ref{prop: general distance with thetas}. It should be compared with Proposition \ref{clry: FM-norm Dirac min prob} for the case $\mu\in\CP(S)$:
\begin{corollary}\label{clry:FM dist Dirac pos measure}
	Let $x\in S$ and $\mu\in\CM^+(S)$. Then
	\[
	\bigl\|\delta_x-\mu\bigr\|^*_\FM = \sup_{\theta\in[-1,1]} \bigl\langle \delta_x-\mu,(-\ind)\vee(\theta- d(x,\cdot))\bigr\rangle.
	\]
\end{corollary}
\noindent Proposition \ref{clry: FM-norm Dirac min prob} states that for $\mu\in\CP(S)$ the supremum above is attained at the value $\theta=1$. Such a stronger result can be obtained for general positive measures too.

Let $B(x,r) :=\{ y\in S: d(x,y)<r \}$ be the open ball in $(S,d)$ of radius $r$, centred at $x$. In the following result the function $(x,r)\mapsto \mu(B(x,r))$ plays a key role. It `measures' in a way the mass distribution of $\mu$ over space.
\begin{theo}\label{prop:FM-norm Dirac minus positive}
	Let $x\in S$ and $\mu\in\CM^+(S)$. Then
	\[
	\bigl\| \delta_x-\mu\bigr\|_\FM^* = \bigl\langle \delta_x-\mu, (-\ind)\vee (\theta_0 - d(x,\cdot)) \bigr\rangle,
	\]	
	where $\theta_0=\theta_0(x) := \bigl(2\wedge\inf\bigl\{r\geq 0: \mu\bigl(B(x,r)\bigr)\geq 1\bigr\}\bigr) - 1\in  [-1,1]$.
\end{theo}
\noindent (We use the convention, that $\inf\emptyset = +\infty$).
\begin{proof}
	Define $\phi(\theta) := \langle\delta_x-\mu, (-\ind)\vee(\theta-d(x,\cdot))\rangle$. Note that $\phi$ is continuous. In view of Corollary \ref{clry:FM dist Dirac pos measure} we have to maximize $\phi$ over $[-1,1]$. Let $\theta,\tilde{\theta}\in[-1,1]$, such that $\theta>\tilde{\theta}$. Then
	\begin{align*}
		\phi(\theta) - \phi(\tilde{\theta}) & = \theta\ -\  \int_{B(x,\theta+1)} \theta-d(x,y)\,d\mu(y)\ -\ \int_{S\setminus B(x,\theta+1)} -\ind\,d\mu\\
		& \qquad -\tilde{\theta}\ +\ \int_{B(x,\tilde{\theta}+1)} \tilde{\theta}-d(x,y)\,d\mu(y)\ +\ \int_{S\setminus B(x,\tilde{\theta}+1)} -\ind\,d\mu\\
		& = (\theta-\tilde{\theta})\bigl( 1- \mu\bigl(B(x,\tilde{\theta}+1)\bigr) \bigr)\ -\ \int_{B(x,\theta+1)\setminus B(x,\tilde{\theta}+1)} \theta - d(x,y) + 1\, d\mu(y).
	\end{align*}
	Here we used that
	\begin{align}
	\int_{S\setminus B(x,\theta+1)} \ind\,d\mu - \int_{S\setminus B(x,\tilde{\theta}+1)} \ind\,d\mu\ &=\ - \bigl[ \mu(B(x,\theta+1)\bigr) - \mu\bigl( B(x,\tilde{\theta}+1)\bigr) \bigr]\\
	& =\ - \mu(B(x,\theta+1)\setminus B(x,\tilde{\theta}+1)\bigr).
	\end{align}
	Therefore, $\phi(\theta)>\phi(\tilde{\theta})$ if and only if
	\begin{equation}\label{eq:condition increasing phi}
		(\theta-\tilde{\theta})\bigl(1- \mu\bigl(B(x,\tilde{\theta}+1)\bigr)\bigr)\ >\  \int_{B(x,\theta+1)\setminus B(x,\tilde{\theta}+1)} \theta - d(x,y) + 1\, d\mu(y).
	\end{equation}
	The function $\tilde{\theta}\mapsto \mu(B(x,\tilde{\theta}+1))$ is non-decreasing, since $\mu$ is a positive measure. According to the definition of $\theta_0$, $\mu(B(x,\tilde{\theta}+1))\geq 1$ for all $\tilde{\theta}>\theta_0$. In that case, inequality \eqref{eq:condition increasing phi} cannot hold, because the right-hand side is non-negative. We conclude that $\phi$ is non-increasing on $(\theta_0,1]$. (If $\theta_0=1$, then $(\theta_0,1]=\emptyset$ and this statement is true trivially.)
	
	We claim that $\phi$ is strictly increasing on $[-1,\theta_0)$.
	If $\theta_0=-1$, then $[-1,\theta_0)=\emptyset$ and there is nothing to prove. So assume $\theta_0>-1$.
	To prove the claim in this case, take $\theta,\tilde{\theta}\in [-1,\theta_0)$, $\theta>\tilde{\theta}$. For all $y\in B(x,\theta+1)\setminus B(x,\tilde{\theta}+1)$ one has $d(x,y)\geq \tilde{\theta}+1$, so $\theta-\tilde{\theta} \geq \theta -d(x,y) +1$.  So if the condition
	\begin{equation}\label{eq:condition 2}
		(\theta-\tilde{\theta})\bigl( 1- \mu\bigl(B(x,\tilde{\theta}+1)\bigr) \bigr)\ >\ \int_{B(x,\theta+1)\setminus B(x,\tilde{\theta}+1)} \theta - \tilde{\theta}\, d\mu(y)
	\end{equation}
	is satisfied, then also condition \eqref{eq:condition increasing phi}. Since $\theta>\tilde{\theta}$, condition \eqref{eq:condition 2} holds if and only if
	\[
	1- \mu\bigl(B(x,\tilde{\theta}+1)\bigr) > \mu\bigl( B(x,\theta+1)\setminus B(x,\tilde{\theta}+1) \bigr),
	\]
	which is equivalent to the condition
	\begin{equation}\label{eq:condition 3}
		\mu\bigl(B(x,\theta+1)\bigr) < 1.
	\end{equation}
	Thus, if condition \eqref{eq:condition 3} holds, then \eqref{eq:condition increasing phi} is satisfied and $\phi(\theta)> \phi(\tilde{\theta})$. By definition of $\theta_0$, \eqref{eq:condition 3} holds for all $\theta<\theta_0$. Thus, $\phi$ is strictly increasing on $[-1, \theta_0)$.
	
	Because $\phi$ is continuous on $[-1,1]$, strictly increasing on $[-1,\theta_0)$ and non-increasing on $(\theta_0,1]$, $\phi$ attains its maximum value at $\theta_0$.
\end{proof}

\begin{remark}
	Paradoxically, in the special case that $\mu\in\CP(S)$, Proposition \ref{clry: FM-norm Dirac min prob} claims that the value of the norm $\|\delta_x-\mu\|_\FM^*$ can be obtained by taking $\theta=1$, instead of $\theta=\theta_0$, which is possibly less than $1$. However, the proof of Proposition \ref{prop:FM-norm Dirac minus positive} shows that the {\it minimum} value for $\theta$ at which the maximum of the function $\phi$ is attained, which equals the stated Fortet-Mourier norm, is $\theta_0$. Since $\phi$ is non-increasing on $(\theta_0,1]$, there must exist also a maximum value at which this (same) maximum value is attained, say $\theta_1$. In case $\mu$ is a probability measure, Proposition \ref{clry: FM-norm Dirac min prob} shows that $\theta_1=1$. So there is no contradiction between the result of Theorem \ref{prop:FM-norm Dirac minus positive} and Proposition \ref{clry: FM-norm Dirac min prob}.
\end{remark}

Let us collect some immediate consequences of Theorem \ref{prop:FM-norm Dirac minus positive}.
\begin{corollary}
	Let $x\in S$ and $\mu\in\CM^+(S)$ with $\mu(S)<1$. Then $\theta_0=1$ and 
	\[
	\bigl\| \delta_x-\mu\bigr\|_\FM^* = \bigl\langle \delta_x-\mu, (-\ind)\vee (1 - d(x,\cdot)) \bigr\rangle = 1-\mu(S) + \bigl\langle \mu, 2\wedge d(x,\cdot)\bigr\rangle.
	\]
\end{corollary}

\begin{corollary}\label{clry:FM-distance wiighted Dirac pos meas}
	Let $x\in S$, $\alpha>0$ and $\mu\in\CM^+(S)$. Then 
	\[
	\bigl\| \alpha\delta_x-\mu\bigr\|_\FM^* = \bigl\langle \alpha\delta_x-\mu, (-\ind)\vee (\theta^\alpha_0 - d(x,\cdot)) \bigr\rangle = \alpha\theta^\alpha_0 - \bigl\langle \mu, (-\ind)\vee (\theta^\alpha_0-d(x,\cdot)) \bigr\rangle,
	\]
	where $\theta^\alpha_0 = \theta^\alpha_0(x) = 2\wedge\inf\bigl\{ r\geq 0: \mu\bigl(B(x,r)\bigr) \geq\alpha \} - 1\in[-1,1]$.
\end{corollary}
\begin{proof}
	One has $\| \alpha\delta_x-\mu\|_\FM^* = \alpha\|\delta_x - \alpha^{-1}\mu\|_\FM^*$. Now apply Theorem \ref{prop:FM-norm Dirac minus positive} to the measure $\alpha^{-1}\mu$ instead of $\mu$.	
\end{proof}

The following corollary nicely generalizes the expression for the well-known Fortet-Mourier distance between two Dirac measures, see \eqref{eq:well-known norms}:
\begin{corollary}\label{clry:FM-distance two weighted single Diracs}
	Let $x,y\in S$ and $\alpha, \beta>0$. Then
	\begin{equation}\label{eq:FM-distance weighted Diarcs}
		\bigl\| \alpha\delta_x - \beta\delta_y \bigr\|_\FM^* = |\alpha-\beta| + (\alpha\wedge \beta)(2\wedge d(x,y)).
	\end{equation}
\end{corollary}
\begin{proof}
	Corollary \ref{clry:FM-distance wiighted Dirac pos meas} yields
	\begin{equation}\label{eq:first expression weighted Diracs}
		\bigl\| \alpha\delta_x - \beta\delta_y \bigr\|_\FM^* = \alpha \theta^\alpha_0 - \beta \bigl( (-1)\vee (\theta^\alpha_0 - d(x,y)) \bigr),
	\end{equation}
	with $	\theta^\alpha_0 = \bigl(2\wedge \inf\bigl\{ r\geq0: \delta_y\bigl( B(x,r)\bigr) \geq \alpha/\beta \bigr\}\bigr)-1\in[-1,1]$. One easily checks that
	\[
	\bigl\{ r\geq 0: \delta_y\bigl( B(x,r)\bigr) \geq \alpha/\beta \bigr\} = \begin{cases} 
		\bigl\{ r:  r>d(x,y) \bigr\}, & \mbox{if}\ \alpha/\beta \leq 1,\\
		\emptyset,  & \mbox{if}\ \alpha/\beta >1.
	\end{cases}
	\]
	Therefore, 
	\[
	\theta^\alpha_0 = \begin{cases} 
		d(x,y) - 1, & \mbox{if}\ \alpha\leq \beta\ \mbox{and}\ d(x,y)<2,\\
		1, & \mbox{otherwise}. \end{cases}
	\]
	From \eqref{eq:first expression weighted Diracs} we find that \begin{align}
		\bigl\|\alpha\delta_x -\beta\delta_y\bigr\|_\FM^* & = {\begin{cases}
				\alpha(d(x,y)-1) + \beta, & \mbox{if}\ \alpha\leq\beta\ \mbox{and}\ d(x,y)<2,\\
				\alpha - \beta ((-1)\vee 1-d(x,y)), & \mbox{otherwise},
		\end{cases}}\nonumber\\
		& = \begin{cases}
			\beta - \alpha + \alpha d(x,y), & \mbox{if}\ \alpha\leq\beta\ \mbox{and}\ d(x,y)<2,\\
			\alpha - \beta\bigl( 1- 2\wedge d(x,y)\bigr), & \mbox{otherwise},
		\end{cases}\nonumber\\
		& = \begin{cases}
			|\alpha - \beta| + \alpha d(x,y), & \mbox{if}\ \alpha\leq\beta\ \mbox{and}\ d(x,y)<2,\\
			\alpha - \beta  + \beta(2\wedge d(x,y)), & \mbox{otherwise},
		\end{cases}\nonumber\\
		& = |\alpha-\beta| + (\alpha\wedge\beta)(2\wedge d(x,y)).\label{eq:final expression}
	\end{align}
	Here we used \eqref{eq:rewriting function} in the second step. To get to \eqref{eq:final expression}, for the case $\alpha\leq \beta$ and $d(x,y)\geq 2$, we used that $2(\alpha\wedge\beta) = \alpha+\beta - |\alpha-\beta|$.
\end{proof}

\begin{remark}
	Without the use of the results that we presented, one could estimate as follows. Put $z:=x$ if $\alpha\geq \beta$ and $z:=y$ if $\alpha<\beta$. Then
	\begin{align*}
		\bigl\|\alpha\delta_x -\beta\delta_y\bigr\|_\FM^* & =\bigl\| (\alpha\wedge\beta)(\delta_x-\delta_y)\ + \ |\alpha-\beta| \delta_z\bigr\|_\FM^*\\
		& \leq (\alpha\wedge\beta) \bigl\| \delta_x-\delta_y\bigl \|_FM^* + |\alpha-\beta| \|\delta_z\|_\FM^*\\
		& = (\alpha\wedge\beta)(2\wedge d(x,y) + |\alpha-\beta|.
	\end{align*}
	The point of Corollary \ref{clry:FM-distance two weighted single Diracs} is, that equality holds.
\end{remark}

We could not obtain an explicit expression for the suprema in \eqref{eq:FM-norm via BP} or \eqref{eq:FM-norm via thetas}, like \eqref{eq:FM-distance weighted Diarcs}, when $\nu$ is a weighted sum of two or more Dirac measures. The distance can be computed though in those cases, by algorithms that we shall exhibit in the next section.

\section{Distance to positive molecular measures -- an algorithmic approach}
\label{sec:algorithms}

We are now concerned with computing $\|\nu-\mu\|_\FM^*$ where $\nu\in\Mol^+(S)$ and $\mu\in\CM^+(S)$. Explicit expressions, like those presented for $\nu=\alpha\delta_x$ in the previous section, could not be obtained. It is possible to compute the norm in particular cases, most importantly when $\mu$ is also a positive molecular measure. Put otherwise, we shall provide an exact algorithm to compute $\|\mu\|_\FM^*$ for any $\mu\in\Mol(S)$. Note we assume the generality of $(S,d)$ being a metric space. Thus, our results provide a substantial generalisation of both \cite{Jablonski-MCz:2013} and \cite{Sriperumbudur_ea:2012}.
\vskip 0.1cm

Let $\nu\in\Mol^+(S)$ and $\mu\in\CM^+(S)$ and write $\nu=\sum_{i=1}^N \alpha_i\delta_{x_i}$ with all $x_i$ distinct and $\alpha_i>0$. Put $P:= \supp(\nu)=\{x_1,\dots,x_N\}$ and view $P$ as a metric space for the restriction of $d$ to $P$. Define $\iota:B^P_\FM\to[-1,1]^N$ by $\iota(f):= (f(x_1),\dots,f(x_N))$. Since for any $f,g\in B^P_\FM$ and $j\in\{1,\dots,N\}$,
\[
\bigl| \iota(f)_j - \iota(g)_j\bigr| \leq \|f-g\|_\infty \leq \|f-g\|_\FM,
\]
$\iota$ is a non-expansive map when $[-1,1]^N$ is equipped with the $\max$-distance or the Euclidean metric. For any $\tau\in\CM(S)$, define
\begin{equation}
	\psi_\tau: [-1,1]^N \to \RR: \theta \mapsto \bigl\langle \tau, (-\ind)\vee\bigvee_{i=1}^N (\theta_i - d(x_i,\cdot)) \bigr\rangle.
\end{equation}
The motivation for studying this function is given by Theorem \ref{prop: general distance with thetas}:
\begin{equation}\label{eq:key equation maximization}
	\bigl\|\nu-\mu\bigr\|_\FM^* = \sup_{\theta\in [-1,1]^N} \psi_{\nu-\mu}(\theta) = \sup_{f\in B^P_\FM} \psi_{\nu-\mu}\bigl(\iota(f)\bigr).
\end{equation}
\begin{prop}\label{prop: props psi tau}
	If $[-1,1]^N$ is equipped with the $\max$-distance or the Euclidean metric, then $\psi_\tau$ is Lipschitz continuous and $|\psi_\tau|_L\leq \|\tau\|_\TV$.
	For $\tau\in\CM^+(S)$, $\psi_\tau$ is convex on $[-1,1]^N$. 
\end{prop}
\begin{proof}
	Let $\theta,\tilde{\theta}\in[-1,1]^N$. Recall the definition of $h_\theta$ in \eqref{def:h_theta}. According to Lemma \ref{lem:properties max operator} {\it (ii)}, for every $x\in S$ one has
	\begin{equation}
		\bigl| h_\theta(x)-h_{\tilde{\theta}}(x) \bigr|\ \leq\  \max_{1\leq i\leq n} |\,\theta_i -\tilde{\theta}_i|.
	\end{equation}
	This yields that for any $\tau\in\CM(S)$,
	\begin{align}
		\bigl| \psi_\tau(\theta) - \psi_\tau(\tilde{\theta}) \bigr| &\ =\ \bigl| \bigl\langle \tau, h_\theta-h_{\tilde{\theta}}\bigr\rangle \bigr|\ \leq\ 
		\|\tau\|_\TV \|h_\theta- h_{\tilde{\theta}}\|_\infty\nonumber\\
		&\  \leq \  \|\tau\|_\TV \max_{1\leq i\leq n} |\,\theta_i -\tilde{\theta}_i|\ \leq\ \|\tau\|_\TV\; \left(\sum_{i=1}^N |\theta_i-\tilde{\theta}_i|^2\right)^{1/2}.
	\end{align}
	So, $\psi_\tau$ is Lipschitz with $|\psi_\tau|_L\leq \|\tau\|_\TV$.
	
	Now assume that $\tau\in\CM^+(S)$. Let $\theta,\tilde{\theta}\in[-1,1]^N$ and $0\leq t\leq 1$. Then, using Lemma \ref{lem:properties max operator} {\it(ii)} and the positivity of $\tau$ to get to inequalities \eqref{eq:inequality max 1} and \eqref{eq:inequality max 2} below, we arrive at
	\begin{align}
		\psi_\tau(t\theta + (1-t)\tilde{\theta})\ &=\ 
		\bigl\langle \tau, (-\ind)\vee \bigvee_{i=1}^N \bigl(t(\theta_i - d(x_i,\cdot))+ (1-t)(\tilde{\theta}_i - d(x_i,\cdot))\bigr) \;\bigr\rangle\nonumber\\
		&\leq\ \bigl\langle \tau, \bigl( t(-\ind)+ (1-t)(-\ind)\bigr)\vee \left( 
		t\bigvee_{i=1}^N (\theta_i - d(x_i,\cdot)) + (1-t) \bigvee_{i=1}^N (\tilde{\theta}_i - d(x_i,\cdot))
		\right)\;\bigr\rangle\label{eq:inequality max 1}\\
		&\leq\ \bigl\langle \tau, t\left( (-\ind)\vee \bigvee_{i=1}^N (\theta_i - d(x_i,\cdot))  \right) + (1-t)\left( (-\ind)\vee \bigvee_{i=1}^N(\tilde{\theta}_i - d(x_i,\cdot))  \right)\;\bigr\rangle\label{eq:inequality max 2}\\
		&=\ t\psi_\tau(\theta)\ +\ (1-t)\psi_\tau(\tilde{\theta}).\nonumber
	\end{align}
	Thus, $\psi_\tau$ is convex on $[-1,1]^N$. 
\end{proof}
\noindent Because $\psi_{\nu-\mu}$ is continuous, the suprema in \eqref{eq:key equation maximization} are attained on the compact sets $[-1,1]^N$ and $\iota(B^P_\FM)$, respectively. 
\vskip 0.1cm

For general signed measure $\tau$, $\psi_\tau$ is the difference of the convex functions $\psi_{\tau^+}$ and $\psi_{\tau^-}$. Consequently, no particular `convexity properties' of $\psi_\tau$ can be claimed. However, for $\tau=\nu-\mu$ with $\nu$ and $\mu$ positive measures as above, one can derive:
\begin{prop}\label{prop: concavity of psi on ball}
	Let $\nu\in\Mol^+(S)$ with $P:=\supp(\nu)=\{x_1,\dots, x_N\}$ and $\mu\in\CM^+(S)$. Then $\psi_{\nu}$ is `linear' on $\iota(B^P_\FM)\subset [-1,1]^N$. In particular, $\psi_{\nu-\mu}=\psi_\nu-\psi_\mu$ is concave.
\end{prop}
\begin{proof}
	One has $\psi_{\nu-\mu}=\psi_\nu-\psi_\mu$. $\psi_\mu$ is convex on $[-1,1]^N$, according to Proposition \ref{prop: props psi tau}, so $-\psi_\mu$ is concave.
	We conclude by showing that $\psi_\nu$ is concave on $\iota(B^P_\FM)$. For general $\theta\in[-1,1]^N$ one has
	$h_\theta(x_i)\geq \theta_i$ for all $i\in\{1,\dots,N\}$. However, for $\theta\in \iota(B^P_\FM)$ there is $f\in B^P_\FM$ such that $\theta=\iota(f)$.  Therefore, according to Lemma \ref{lem:bounding from below by extreme point} {\it (ii)}, $h_\theta(x_i)=f(x_i)$ for all $i$. Thus, for any $f\in B^P_\FM$, writing $f_i=\iota(f)_i=f(x_i)$,
	\begin{equation}\label{eq:shape psi pos mol measure}
		\psi_\nu(\iota(f))\ =\ \bigr\langle \sum_{i=1}^N \alpha_i\delta_{x_i}, h_{\iota(f)} \bigr\rangle\ =\ \sum_{i=1}^N \alpha_i f_i.
	\end{equation}
	So $\psi_\nu$ on $\iota(B^P_\FM)$ is the restriction of a linear functional on $\RR^N$ to the convex subset $\iota(B^P_\FM)$. In particular, $\psi_\nu$ is concave on $\iota(B^P_\FM)$.	
\end{proof}

Thus, in view of \eqref{eq:key equation maximization} and Proposition \ref{prop: concavity of psi on ball} the problem of computing the Fortet-Mourier distance $\|\nu-\mu\|_\FM^*$ for $\nu\in\Mol^+(S)$ and $\mu\in \CM^+(S)$ is equivalent to the problem of maximizing the concave and Lipschitzian function $\psi_{\nu-\mu}$ over the compact convex set $\iota(B^P_\FM)$ in $\RR^N$, where $P=\supp(\nu)\subset S$, consisting of $N$ distinct points. This is again equivalent to minimizing the convex function $-\psi_{\nu-\mu}$ over $\iota(B^P_\FM)$.

Minimization of convex functions has been widely studied (cf. eg. \cite{Bertsekas:2015,Boyd-VandenBerghe:2004,Rockafellar:1972}) and a wide variety of algorithms have been developed for convex minimization in the field of convex optimization. These problems can be solved highly efficiently by now. Boyd and VandenBerghe even state (\cite{Boyd-VandenBerghe:2004} p.8):
\vskip 0.1cm

\noindent `{\it With only a bit of exaggeration, we can say that, if you formulate a practical problem as a convex optimization problem, then you have solved the original problem.}'
\vskip 0.1cm

\noindent Therefore, we consider the theoretical side of computing distances of the form $\|\nu-\mu\|_\FM^*$ with $\nu\in\Mol^+(S)$ and $\mu\in\CM^+(S)$ as solved. 

Of course, in a practical setting, the implementation of the convex optimization algorithm of choice for a specific measure $\mu$ may require additional practical issues to be resolved. For example, one must be able to compute $\psi_\mu$ (approximately). In the following section we shall consider the special case where $\mu\in\Mol^+(S)$, i.e. computing the Fortet-Mournier norm of a molecular measure. But let us provide another example first.

\begin{example}
	Let $S=[0,\infty)$, equipped with the Euclidean metric and take \[
	\mu := e^{-x^2}dx, \qquad \nu := \smfrac{1}{2}\delta_0 + 2\delta_{\smfrac{1}{3}} + \smfrac{1}{5}\delta_{\smfrac{1}{2}} + \smfrac{1}{3}\delta_3.
	\]
	We implemented an algorithm to mimimize $-\psi_{\nu-\mu}$ over $\iota(B^P_\FM)$, with $P=\bigl\{0, \smfrac{1}{3},\smfrac{1}{2}, 3\bigr\}$, using the MATLAB `\texttt{fmincon}' function, see Appendix \ref{app: FM dual norm fmincon}. It resulted in
	\[
	\bigl\|\nu-\mu\bigr\|_\FM^* = -\bigl(-\psi_{\nu-\mu}(f)\bigr) \approx 2.3921\dots, \qquad\mbox{with}\ f = \bigl[ 1,\ 1,\ \smfrac{5}{6},\ 1\bigr].
	\]
\end{example}

\section{Computing the Fortet-Mourier norm of a molecular measure}
\label{sec:computing FM-norm molecular}

We conclude by specializing to the particular case where both $\nu$ and $\mu$ are positive molecular measures. That is, we show how to compute $\|\tau\|_\FM^*$ for $\tau\in\Mol(S)$. We present two ways to proceed: one by specializing the results of the previous section and one that is special to this specific case. Each has its benefits and drawbacks, which we shall discuss. We start with the latter method. We note that \cite{Jablonski-MCz:2013} (see also \cite{Evers:thesis} Appendix) provides an algorithm to compute $\|\tau\|_\FM^*$ when $S=\RR$ or an interval therein. \cite{Sriperumbudur_ea:2012} exhibits a method that works for general space $S$, but $\tau$ must be the difference of two empirical measures. That is, the coefficients of the Diracs are quite specific. Before starting any further considerations, note that if $\tau$ or $-\tau$ is positive, then $\|\tau\|_\FM^* =\|\tau\|_\TV = \sum_i |\alpha_i|$ if $\tau =\sum_i\alpha_i\delta_{x_i}$. Thus, we shall assume $\tau^+\neq0$ and $\tau^-\neq 0$.
\vskip 0.1cm

\begin{figure}[h]
	\begin{tabular}{ll}
		\includegraphics[trim=0 0 3.2cm 0, clip,scale=0.5]{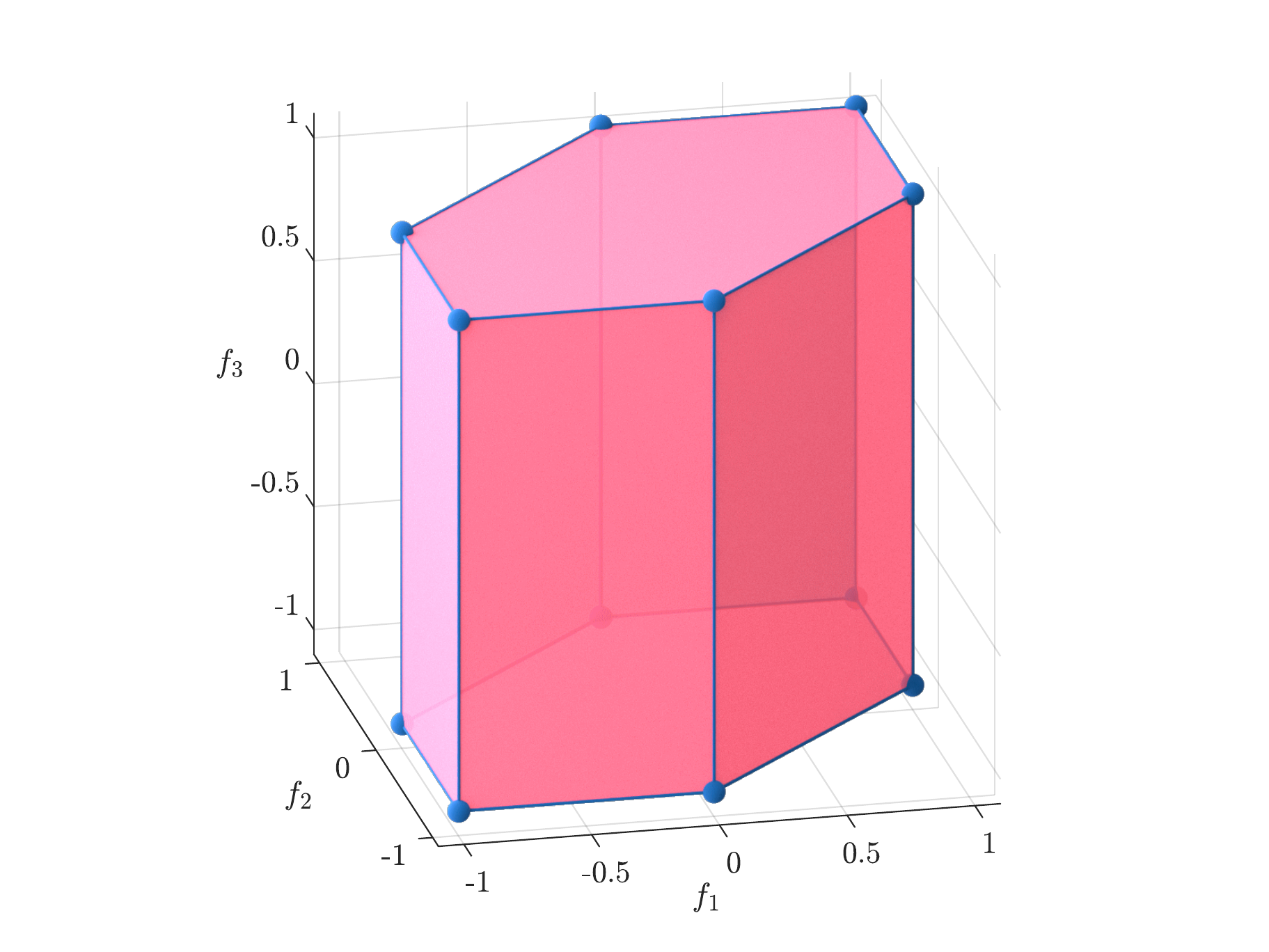}	& 
		\includegraphics[trim=2.2cm 0 0 0, clip, scale=0.5]{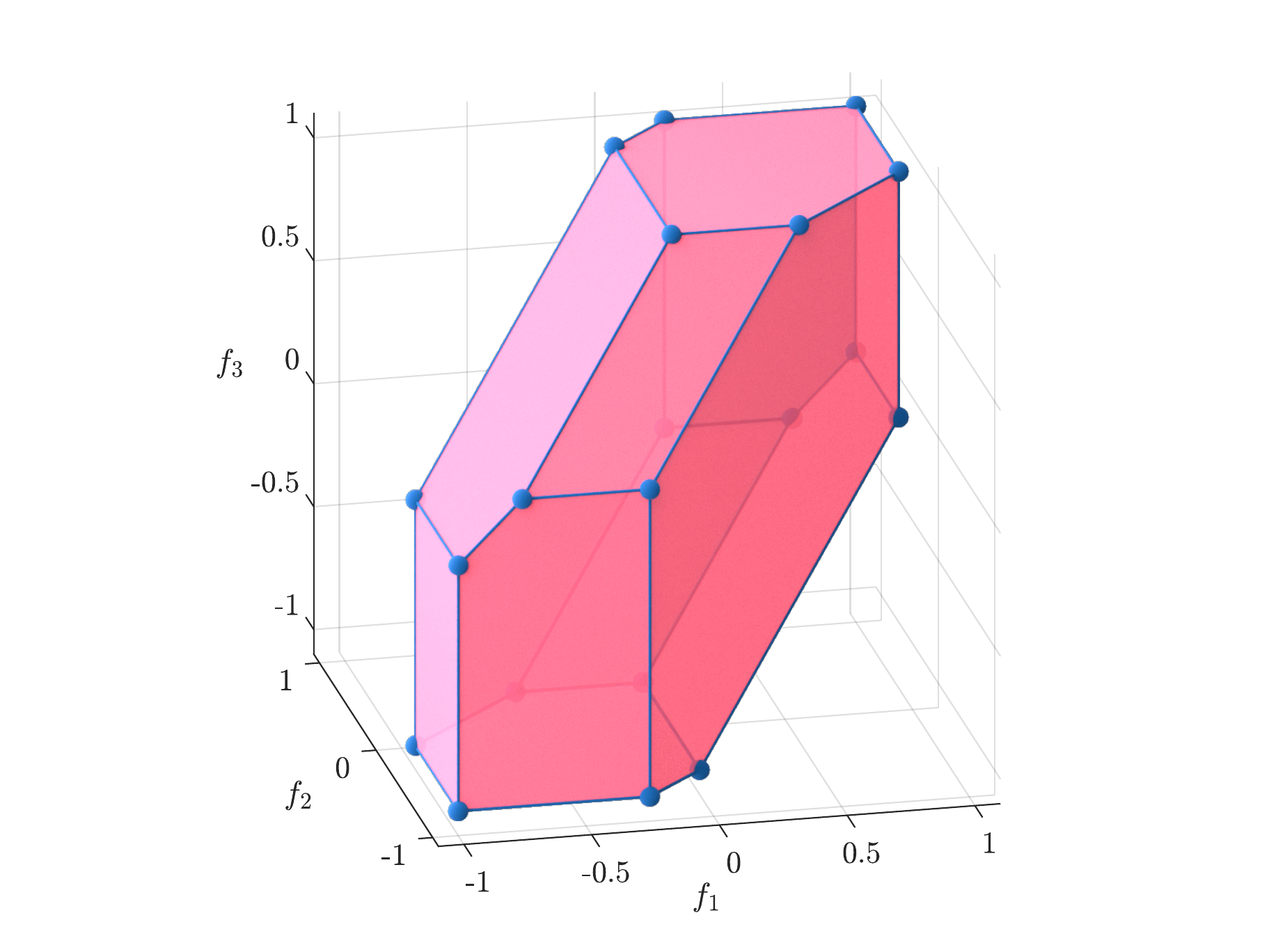}\\
	\end{tabular}
	\caption{{\it The unit ball $B^P_\FM$ for the norm $\|\cdot\|_\FM$ on the space $\BL(P,d)$, where $P=\{x_1,x_2,x_3\}$ and $f\in\BL(P,d)$ is represented by $(f_i)\in\RR^3$ with $f_i:=f(x_i)$. The defining conditions are given by \eqref{eq:standard lin programming representation}. The metric $d$ is defined by $d_{ij}:=d(x_i,x_j)$ and differs for the two cases shown.\\
			Left: $d_{12}=1$, $d_{13}=2$, $d_{23}=3$. Right: $d_{12}=0.75$, $d_{13}=1$, $d_{23}=1.25$.}\label{fig:FM-balls}}	
\end{figure}

As before, let $(S,d$) be a metric space and $P=\{x_1,\dots x_n\}$ a set of $n$ distinct points in $S$. $P$ inherits the metric structure of $S$, by restriction. Put $d_{ij}:=d(x_i,x_j)$. It is readily verified that
\begin{align}
	\iota(B^P_\FM)\ &=\ \bigcap_{1\leq k\leq n} \bigl\{ f\in\RR^n: |f_k|\leq 1 \bigr\}\ \ \cap\ \bigcap_{1\leq i<j\leq n} \bigl\{ f\in\RR^n: |f_i-f_j|d_{ij}^{-1}\leq 1 \bigr\}\nonumber\\
	&=\ \bigcap_{1\leq k\leq n} \bigl\{ f_k\leq1 \bigr\}\cap \bigl\{ -f_k\leq 1\bigr\}\ \ \cap\ \bigcap_{1\leq i<j\leq n} \bigl\{ (f_i-f_j)d_{ij}^{-1} \leq 1\bigr\} \cap \bigl\{ (f_j-f_i)d_{ij}^{-1} \leq 1\bigr\}.\label{eq:standard lin programming representation}
\end{align}
Expression \eqref{eq:standard lin programming representation} is in the form of the standard linear programming representation of the domain of the objective function as an intersection of finitely many half-spaces (see e.g. \cite{Bertsimas:1997}). Figure \ref{fig:FM-balls} shows two unit balls $B^P_\FM$, with $P$ consisting of three points, for two different metrics.

%%%%%%%%%%%%%%%%%%%%%%%%%%%%%%%%%
%%%%%%%%%%%%%%%%%%%%%%%%%%%%%%%%%

Write $\tau = \sum_{i=1}^n \alpha_i \delta_{x_i}$ with $0\neq\alpha_i\in \RR$ and $x_i\in S$, and put $P:=\supp(\tau) = \{ x_1,\dots x_n\}$. Lemma \ref{lem:bounding from below by extreme point} {\it (ii)} implies that restriction to $P$ gives a surjective map from $B^S_\FM$ onto $B^P_\FM$. Therefore,
\begin{equation}\label{eq:translation FM-norm lin prog problem}
	\| \tau \|_\FM^*\ =\ \sup_{g\in B^S_\FM} \langle \tau, g,\rangle \ =\ \sup_{f\in B^P_\FM} \langle \tau|_P, f\rangle = \max_{f\in \iota(B^P_\FM)}\ \sum_{i=1}^n \alpha_i f_i.
\end{equation}

Thus, $\|\tau\|_\FM^*$ can be computed using one of the many existing -- very efficient -- optimization algorithms that use linear programming, such as Gurobi and CPLEX, or the built-in `\texttt{linprog}' function in MATLAB, using the standard domain description \eqref{eq:standard lin programming representation}. In these algorithms there is an initial step in which an extreme point of the domain is sought to start the search for the optimum. Here, $\pm(1,\dots,1)$ are always extreme points. If $\sum_i\alpha_i\geq 0$ one may start at $(1,\dots, 1)$. If $\sum_i\alpha_i< 0$,one may start at $(-1,\dots, -1)$. This reduces the number of vertices of $B^P_\FM$ that needs to be examined by the optimization algorithm in the worst case by a factor two.
%\newpage

\begin{example}
	Let $P=\{x_1, x_2, x_3\}$ and $d_{12}=1$, $d_{13}=2$ and $d_{23}=3$. Notice that the prescribed distances are consistent with the triangle inequality. Take $\alpha_1=1$, $\alpha_2=-\smfrac{1}{3}$ and $\alpha_{3} = -\smfrac{2}{3}$. Equation \eqref{eq:FM-norm 1-2 Diracs} gives an explicit result in this case:
	\[
	\bigl\| \sum_i \alpha_i\delta_{x_i} \bigr\|_\FM^*\ =\ |\alpha_2| (2\wedge d_{12}) + (1-|\alpha_2|)(2\wedge d_{13}) \ =\ \frac{1}{3}\cdot 1 + \frac{2}{3}\cdot 2 = \frac{5}{3}.
	\]
	We implemented an algorithm for computing the Fortet-Mourier norms of molecular measures in MATLAB, using the `\texttt{linprog}' algorithm, see Appendix \ref{app:FM dual norm}. It returned the same result for the norm, but also an $f$ at which the optimum is attained. In this case $f= (1, 0, -1)$. This $f$ corresponds precisely to the function $(-\ind)\vee(1-d(x_1,\cdot))$ that appears in the theoretical result, Proposition \ref{clry: FM-norm Dirac min prob}.
\end{example}

The above linear programming algorithm for computing $\|\tau\|_\FM^*$ has as domain for the objective function a polygon in $\RR^n$, where $n$ is the number of points in the support of $\tau\in\Mol(S)$. The dimensionality of the optimization problem can be reduced by halve, by resorting to the results of Section \ref{sec:algorithms}. The number of points in the support of either $\tau^+$ or $\tau^-$ is less than $n/2$ or both have precisely $n/2$ points in their support. The one with the least number of points, say $\tau^-$ with $N\leq n/2$ points, can play the role of $\nu$ in Section \ref{sec:algorithms}, while the other takes up the role of $\mu$, simply because $\|\tau\|_\FM^*=\|-\tau\|^*_\FM$.

Thus, with $\nu=\tau^-$ (say) and $P=\supp(\nu)=\{x_1,\dots, x_N\}$,  
\[
\|\tau\|_\FM^* = \bigl\| \tau^- - \tau^+\bigr\|_\FM^* = \max_{ f\in\iota(B^P_\FM)} \psi_{-\tau}(f) = - \min_{f\in\iota(B^P_\FM)} \bigl( -\psi_{-\tau}(f) \bigr),
\]
according to \eqref{eq:key equation maximization} and the further discussion in Section \ref{sec:algorithms}. The polygonal domain of optimization $\iota(B^P_\FM)$ is still given by \eqref{eq:standard lin programming representation}, but now has reduced dimension $N\leq n/2$. In this setting, $\psi_{\tau^-}$ is `linear' on $\iota(B^P_\FM)$ (Proposition \ref{prop: concavity of psi on ball}). If 
\[
\tau = \sum_{j=1}^{n-N} \beta_j \delta_{y_j}\ -\ \sum_{i=1}^N \alpha_i\delta_{x_i},\qquad \alpha_i, \beta_j>0,
\]
then for $\theta\in\iota(B^P_\FM)$
\begin{equation}
	\psi_{-\tau}(\theta) = \sum_{i=1}^N \alpha_i\theta_i - \sum_{j=1}^{n-N} \beta_j \bigl((-1)\vee \max_{1\leq i\leq N} (\theta_i - d(x_i, y_j)) 
	\bigr).
\end{equation}
So, the reduction in dimension of the domain is to the cost of (part of the) `linearity' of the objective function $-\psi_{-\tau}$ on $\iota(B^P_\FM)$, although it is still convex. Application of one of the existing (efficient) convex optimization algorithms now yields $\|\tau\|_\FM^*$ on a lower dimensional domain.

\section{Computing Dudley norms}
\label{sec:results for Dudly norm}

So far we have been discussing expressions for and the computation of Fortet-Mourier norms only. The dual bounded Lipschitz norm on $\CM(S)$, also known as Dudley norm, or flat metric -- for the associated metric -- is also considered often. It is given by
\[
\|\mu\|_\BL^* := \sup_{f\in B^S_\BL} \langle\mu,f\rangle,\qquad B^S_\BL:=\bigl\{ f\in\BL(S): \|f\|_\infty + |f|_L\leq 1 \bigr\}.
\]

The results presented in Section \ref{sec:dimensional reduction} and Section \ref{sec:FM-distance to single Dirac} do not readily generalize to the $\|\cdot\|_\BL^*$-norm. The main issue is, that the geometric shape of $B^S_\BL$ is more complicated than that of $B^S_\FM$. For example, we could change the $\|\cdot\|_\infty$-norm of a function in $B^S_\FM$ `independently' from its Lipschitz constant. For $B^S_\BL$ that is not possible, since the constraint $\|f\|_\infty + |f|_L\leq 1$ should be maintained.

\begin{figure}[h]
	\includegraphics[scale=0.7]{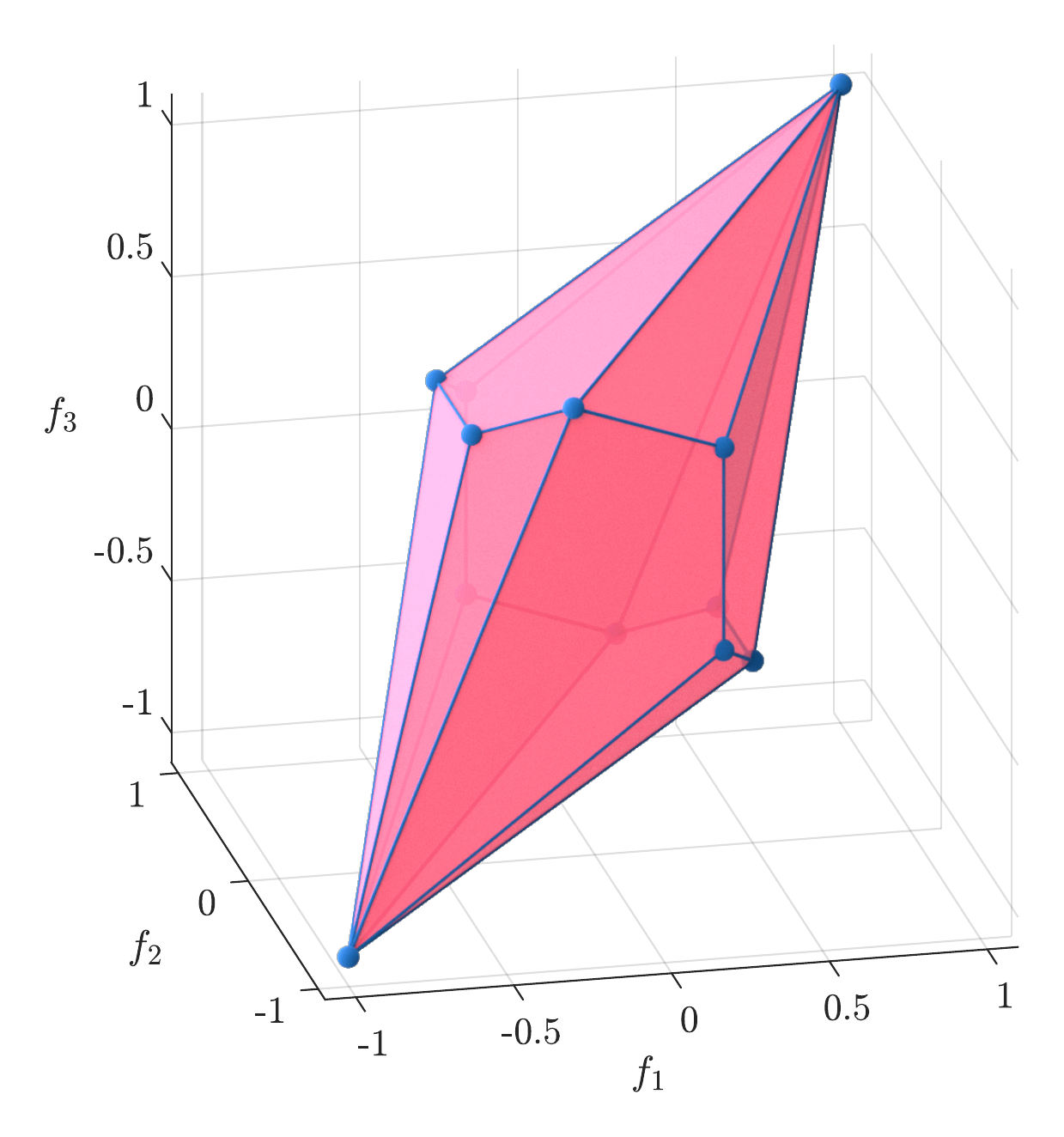}
	\caption{{\it The unit ball $B^P_\BL$ for the norm $\|\cdot\|_\BL$ on the space $\BL(P,d)$, where $P=\{x_1,x_2,x_3\}$ and $f\in\BL(P,d)$ is represented by $(f_i)\in\RR^3$ with $f_i:=f(x_i)$. The defining conditions are given by those in Lemma \ref{lem:description ball BL}, while the metric $d$ is defined by $d_{ij}:=d(x_i,x_j)$ with $d_{12}=1$, $d_{13}=2$, $d_{23}=3$.}\label{fig:BL-ball}}	
\end{figure}

However, the results of Section \ref{sec:computing FM-norm molecular} can be carried over to the $\|\cdot\|_\BL^*$-norm. The corresponding statement of \eqref{eq:translation FM-norm lin prog problem} is still valid with `$\FM$' replaced by `$\BL$'. The description of $\iota(B^P_\BL)$ becomes more awkward though. Without giving the lengthy proof here, we obtained
\begin{lemma}\label{lem:description ball BL}
	Let $n\geq 2$ and $P=\{x_1,\dots,x_n\}\subset S$, consisting of distinct points. Put $d_{ij} := d(x_i,x_j)$. Then
	\begin{align*}
		\iota(B^P_\BL)\ &=\ \bigcap_{1\leq i,j,k\leq n,\ k\neq i\neq j\neq k} \bigl\{ f_k + (f_i-f_j)d_{ij}^{-1}\leq 1\bigr\}\cap \bigl\{ -f_k - (f_i-f_j)d_{ij}^{-1}\leq 1\bigr\}\\
		& \qquad \cap \ \bigcap_{1\leq i, k\leq n,\ i\neq k} 
		\bigl\{ f_k + (f_k-f_i)d_{ik}^{-1}\leq 1\bigr\}\cap \bigl\{ -f_k - (f_k-f_i)d_{ik}^{-1}\leq 1\bigr\}.
	\end{align*}
\end{lemma}
\noindent In Figure \ref{fig:BL-ball} an example is shown of a unit ball $B^P_\BL$ for $P$ consisting of three points. Compare the more complex geometric structure with those of $B^P_\FM$ presented in Figure \ref{fig:FM-balls}.

%%%%%%%%%%%%%%%%%%%%%%%%%%%%%%%%%
%%%%%%%%%%%%%%%%%%%%%%%%%%%%%%%%%

\noindent Thus, one has --  with $P$ as in Lemma \ref{lem:description ball BL} and $0\neq \alpha_i\in \RR$:
\begin{equation}
	\bigl\| \sum_{i=1}^n \alpha_i\delta_{x_i}\bigr\|_\BL^*\ =\  \max_{f\in\iota(B^P_\BL)} \sum_{i=1}^n \alpha_i f_i\ =\ - \min_{f\in \iota(B^P_\BL)} \bigl( - \sum_{i=1}^n \alpha_i f_i \bigr),
\end{equation}
which optimum can again be found by a linear programming algorithm, now by using Lemma \ref{lem:description ball BL} for the standard description of the domain of optimization as intersection of half-spaces.
\vskip 0.2cm

%\newpage
{\bf Acknowledgement.}\ We thank M.A. M\"uller for preparing the graphics showing the unit balls in the spaces of bounded Lipschitz functions over a finite set of points in Figure \ref{fig:FM-balls} and Figure \ref{fig:BL-ball}.

\begin{appendix}
	
	\section{MATLAB implementations}
	
	\subsection{FM-distance between a positive linear combination of Dirac measures and a positive measure}\label{app: FM dual norm fmincon}
	\begin{verbatim}
		% Computes ||\nu-\mu||_FM^* for \nu in Mol^+(S), \mu in M^+(S), 
		% S=(Smin, Smax)\subset\R (possibly Smin=-Inf, Smax=+Inf)
		% \mu=h d\lambda, \mu abs ct wrt Lebesgue measure, \nu=sum a_i delta_{s_i}
		
		Smin = 0;
		Smax = Inf;
		a = [1/2 2 1/5 1/3]; % a=[a_1,...,a_n]
		n = length(a);
		P = [0 1/3 1/2 3]; %=supp(\nu)=[s_1,...,s_n] 
		
		h = @(s) exp(-s^2);  
		integrand = @(s,f)h(s)*max(-1,max(f-abs(P-s)));
		psi = @(f) -dot(a,f)+integral(@(s) integrand(s,f),Smin,Smax,'ArrayValued',true);
		
		A = [];
		for i = 1:n
		    for j = i+1:n
		        B = zeros(1,n);
		        B(i) = abs(P(i)-P(j))^(-1);
		        B(j) = -B(i);
		        A = [A;B];
		    end
		end
		A = [A;eye(n)];
		A = [A;-A];
		
		b = ones(n^2+n,1);
		f0 = zeros(1,n);
		
		[f,val] = fmincon(psi,f0,A,b);
		f
		norm = -val
	\end{verbatim}
	
	\subsection{FM-norm of a linear combination of Dirac  measures}\label{app:FM dual norm}
	\begin{verbatim}
		function [norm,f] = FMdualnorm(a,dist) 
		%a=[a_1 a_2 ... a_n], mu=sum_i=1^n a_i delta_s_i, dist=[d_ij]_i,j=1^n matrix 
		%norm =||mu||_FM^*, f=ext pt for which <mu,f>=norm
		n=length(a);
		A=[];
		for i=1:n
		    for j=i+1:n
		        B=zeros(1,n);
		        B(i)=dist(i,j)^(-1);
		        B(j)=-B(i);
		        A=[A;B];
		    end
		end
		A=[A;eye(n)];
		A=[A;-A];
		b=ones(n^2+n,1); %m=#rows of A=2n+2(n choose 2)=2n+n(n-1)=n^2+n
		minus_a=-a;
		f=linprog(minus_a,A,b);
		norm=a*f;
	\end{verbatim}
	
\end{appendix}

\end{document}